\documentclass{amsart}
\usepackage{amsmath, amsthm, amscd, amsfonts, amssymb, graphicx, color}
\usepackage[bookmarksnumbered]{hyperref}
\newtheorem{thm}{Theorem}[section]

\theoremstyle{definition}
\newtheorem{defn}[thm]{Definition}
\theoremstyle{remark}

\numberwithin{equation}{section}

\begin{document}

\title[Functional equations related to the monomials]{A family of functional equations related to the monomial functions and its stability}%
\author{A. Sousaraei $^{1}$, M. Alimohammady $^{2}$ and A. Sadeghi $^{3}$}%
\address{$^{1}$ Islamic Azad University Branch of Azadshaher, Azadshaher,
Iran.}
\address{$^{2}$ Department of Mathematics, University of Mazandaran,
Babolsar, Iran.}
\address{$^{3}$ Department of Mathematics Tarbiat Modares University Tehran 14115-333 Iran}%
\email{$^{1}$ ali Sousaraei@yahoo.com, $^{2}$ amohsen@umz.ac.ir, $^{3}$ ali.sadeghi@modares.ac.ir}%
\thanks{}%

\begin{abstract}
Our aim of this paper is to study a family of functional equation in
vector and Banach spaces with difference operators, where this
family of functional equation is a general mixed
additive-quadratic-cubic-quartic functional equations. We show that
every function satisfies the our functional equation is a monomial
function with a certain degree. Furthermore, we deal with the
generalized Hyers-Ulam stability of this family of functional
equations in Banach space.
\end{abstract}
\maketitle
\section{Background results}

Throughout this paper, assume that $F=\mathbb{Q}$, $\mathbb{R}$ or
$\mathbb{C}$ and $V$ and $B$ are vector spaces over $F$ and $X$ is a
Banach spaces over $F$. We present some background results based
of paper [1].
A function $a:V\rightarrow B$ is said to be additive provided
$a(x+y)=a(x)+a(y)$ for all $x, y\in V$; in this case it is easily
seen that $a(rx) =ra(x)$ for all $x\in V$ and all $r\in \mathbb{Q}$
.

If $k\in \mathbb{N}$ and $a:V^{k}\rightarrow B$, then we say that a
is $k$-additive provided it is additive in each variable; we say
that a is symmetric provided
$a(x_{1},x_{2},...,x_{k})=a(y_{1},y_{2},...,y_{k})$ whenever
$x_{1},x_{2},...,x_{k}\in V$ and $(y_{1},y_{2},...,y_{k})$ is a
permutation of $(x_{1},x_{2},...,x_{k})$.

If $k\in \mathbb{N}$ and $a:V^{k}\rightarrow B$ is symmetric and
$k$-additive, let $a^{\ast}(x)=a(x,x,...,x)$ for $x\in V$ is and
note that $a^{\ast}(rx)=r^{k}a(x)$ whenever $x\in V$ and $r\in
\mathbb{Q}$. Such a function $a^{\ast}$ will be called a
$\emph{monomial function of degree k}$ (assuming $a^{\ast}\neq 0$).

A function $p:V\rightarrow B$ is called a $\emph{generalized
polynomial}$ (GP) function of degree $m\in \mathbb{N}$ provided
there exist $a_{0}\in B$ and symmetric $k$-additive functions
$a_{k}:V^{k}\rightarrow B$ (for $1\leq k\leq m$) such that
$$p(x)=a_{0}+\sum_{k=1}^{m}a_{k}^{\ast}(x)\ \ \ for\ all\ x\in V,$$
and $a^{\ast}_{m}\neq 0$. In this case
$$p(rx)=a_{0}+\sum_{k=1}^{m}r^{k}a_{k}^{\ast}(x)\ \ \ for\ all\ x\in V\ and\ r\in \mathbb{Q}.$$

Let $B^{V}$ denote the vector space (over $F$) consisting of all
maps from $V$ into $B$. For $h\in V$ define the linear
$\emph{difference}$ operator $\Delta_{h}$ on $B^{V}$ by
$$\Delta_{h}f(x)=f(x+h)-f(x)\ \ \ for\ all\ f\in B^{V} and\ x\in V.$$
Notice that these difference operators commute
($\Delta_{h_{1}}\Delta_{h_{2}}=\Delta_{h_{2}}\Delta_{h_{2}}$ for all
$h_{1},h_{2}\in V$) and if $h\in V$ and $n\in \mathbb{N}$, then
$\Delta_{h}^{n}$--the $n$-th iterate of $\Delta_{h}$--satisfies
$$\Delta_{h}^{n}f(x)=\sum_{k=0}^{n}(-1)^{n-k}(_{k}^{n})f(x+kh)\ \ \ for\ f\in B^{V}\ and\  x,h\in V.$$

The following theorem were proved by Mazur and Orlicz [2] and
[3], and in greater generality by Djokovi$\acute{c}$
[4].
\begin{thm}\label{theo2}
If $n\in \mathbb{N}$ and $f:V\rightarrow B$, then the following are
equivalent.
\begin{enumerate}
\item $\Delta_{h}^{n}f(x)=0$ for all $x, h\in V$.
\item $\Delta_{h_{n}}...\Delta_{h_{1}}f(x)=0$ for all $x,h_{1},...,h_{n}\in V.$
\item $f$ is a GP function of degree at most $n-1$.
\end{enumerate}
\end{thm}
The starting point of the stability theory of functional equations
was the problem formulated by S. M. Ulam in 1940 (see [5]),
during a conference at Wisconsin University:

\emph{Let $(G,.)$ be a group $(B,.,d)$ be a metric group. Does for every $\varepsilon>0$, there exists a $\delta>0$ such that if a
function $f:G\rightarrow B$ satisfies the inequality
$$d(f(xy),f(x)f(y))\leq \delta,\ \ x,y\in G,$$
there exists a homomorphism $g:G\rightarrow B$ such that
$$d(f(x),g(x))\leq\varepsilon,\ \ x\in G?$$}
In 1941, D. H. Hyers [6] gave an affirmative partial answer
to this problem. This is the reason for which today this type of
stability is called Hyers-Ulam stability of functional equation. In
1950, Aoki [7] generalized Hyers' theorem for approximately
additive functions. In 1978, Th. M. Rassias [8] generalized
the theorem of Hyers by considering the stability problem with
unbounded Cauchy differences. Taking this fact into account, the
additive functional equation $f(x+y)=f(x)+f(y)$ is said to have the
Hyers-Ulam-Rassias stability on $(X,Y)$. This terminology is also
applied to the case of other functional equations. For more detailed
definitions of such terminology one can refer to [9] and
[10] Thereafter, the stability problem of functional
equations has been extended in various directions and studied by
several mathematicians
[11-25].

The functional equation
\begin{eqnarray}\label{1}
%
  Q(x+y)+Q(x-y) &=& 2Q(x)+2Q(y)
\end{eqnarray}
is called a quadratic functional equation. In particular, every
solution of the quadratic functional equation is said to be a
quadratic function. The Ulam problem for the quadratic functional
equation was proved by Skof [26,27] for mappings
$f:X\rightarrow Y$, where $X$ is a normed space and $Y$ is a Banach
space. Cholewa [28] noticed that the theorem of Skof is still
true if the relevant domain $X$ is replaced by an Abelian group.
Several functional equations have been investigated by J. M. Rassias
in [29-31], St. Czerwik in [32] and Th. M. Rassias
in [33].

Jun and Kim [34] introduced the following functional equation
\begin{eqnarray}\label{11}
%
f(2x+y)+f(2x-y) &=& 12f(x)+2f(x+y)+2f(x-y)
\end{eqnarray}
and they established the general solution and the generalized
Hyers-Ulam-Rassias stability problem for the functional equation
(1.2). Also a quartic functional equation have been
investigated by J. M. Rassias in [35].

In [36], S. Lee et al. considered the following quartic
functional equation
\begin{eqnarray}\label{n1}
%
  F(2x+y)+F(2x-y)&=&4F(x+y)+4F(x-y)+24F(x)-6F(y)
\end{eqnarray}
and they established the general solution and the stability problem
for the functional equation (1.3). Also the quartic functional
equation have been investigated by J. M. Rassias in [37].

In this paper, we consider a family of functional equations as
follows
\begin{equation}\label{g1}
c_{1}f(ax+by)+c_{2}f(ax-by)=c_{3}f(x)+c_{4}f(y)+c_{5}f(x+y)+c_{6}f(x-y),
\end{equation}
for a given mapping $f:V\rightarrow X$ and for given
$a,b,c_{1},..,c_{6}$ members of $F$, where $a, b, c_{1}$ and $c_{3}$
are nonzero elements and $c_{1}+c_{2}\neq0$. Our aim is to study
functional equation (1.4), where its a general mixed
additive-quadratic-cubic-quartic (1.1), (1.2) and
(1.3) functional equations. We will establish the solution of
the functional equation (1.4) with difference operators. In
fact we proved that $f$ is a monomial function of degree
$\log_{|a|}|\frac{c_{3}+c_{5}+c_{6}}{c_{1}+c_{2}}|$, where its at
most $4$. Furthermore, by using fixed point methods, we are going to
solve the generalized Hyers-Ulam stability problem for the
functional equation (1.4) in Banach space.
\section{Solution of equation (1.4)}

For a given mapping $f:V\rightarrow B$, we define the function
$f:V\times V\rightarrow B$ as follows
$$Df(x,y):=$$
$$c_{1}f(ax+by)+c_{2}f(ax-by)-c_{3}f(x)-c_{4}f(y)-c_{5}f(x+y)-c_{6}f(x-y)$$
for all $x, y\in V.$
\begin{thm}\label{theo3}
Let $f:V\rightarrow B$ be a function such that $Df(x,y)=0$ for all
$x, y\in V$. Then $f$ is a GP function of degree at most $4$.
\end{thm}
\textbf{Proof.} Let $h_{1},h_{2},h_{3},h_{4},h_{5}\in V$ be
arbitrary fixed elements. Letting $x+h_{1}$ and $y-h_{1}$ instead of
$x$ and $y$ in (1.4), we get
$$c_{1}f(ax+by+(a-b)h_{1})+c_{2}f(ax-by+(a+b)h_{1})=$$
$$c_{3}f(x+h_{1})+c_{4}f(y-h_{1})+c_{5}f(x+y)+c_{6}f(x-y+2h_{1})$$
Using this equality and (1.4), we find that
\begin{eqnarray}\label{3}
%
\Delta_{(a-b)h_{1}}c_{1}f(ax+by)+\Delta_{(a+b)h_{1}}c_{2}f(ax-by)&=&
\end{eqnarray}
$$\Delta_{h_{1}}c_{3}f(x)+\Delta_{-h_{1}}c_{4}f(y)+\Delta_{2h_{1}}c_{6}f(x-y).$$
Now, Letting $x+h_{2}$ and $y+h_{2}$ instead of $x$ and $y$ in
(2.1), we get
\begin{eqnarray*}
%
\Delta_{(a-b)h_{1}}c_{1}f(ax+by+(a+b)h_{2})+\Delta_{(a+b)h_{1}}c_{2}f(ax-by+(a-b)h_{2})&=&
\end{eqnarray*}
$$\Delta_{h_{1}}c_{3}f(x+h_{2})+\Delta_{-h_{1}}c_{4}f(y+h_{2})+\Delta_{2h_{1}}c_{6}f(x-y).$$
Using this equality and (2.1), we find that
\begin{eqnarray}\label{4}
%
\Delta_{(a-b)h_{1}}\Delta_{(a+b)h_{2}}c_{1}f(ax+by)+\Delta_{(a+b)h_{1}}\Delta_{(a-b)h_{2}}c_{2}f(ax-by)&=&
\end{eqnarray}
$$\Delta_{h_{1}}\Delta_{h_{2}}c_{3}f(x)+\Delta_{-h_{1}}\Delta_{h_{2}}c_{4}f(y).$$
Letting $x+bh_{3}$ and $y-ah_{3}$ instead of $x$ and $y$ in
(2.1), we get
\begin{eqnarray*}
%
\Delta_{(a-b)h_{1}}\Delta_{(a+b)h_{2}}c_{1}f(ax+by)+\Delta_{(a+b)h_{1}}\Delta_{(a-b)h_{2}}c_{2}f(ax-by+(2ab)h_{3})&=&
\end{eqnarray*}
$$\Delta_{h_{1}}\Delta_{h_{2}}c_{3}f(x+bh_{3})+\Delta_{-h_{1}}\Delta_{h_{2}}c_{4}f(y-ah_{3}).$$
Using this equality and (2.2), we find that
\begin{eqnarray}\label{5}
%
\Delta_{(a+b)h_{1}}\Delta_{(a-b)h_{2}}\Delta_{(2ab)h_{3}}c_{2}f(ax-by)&=&
\end{eqnarray}
$$\Delta_{h_{1}}\Delta_{h_{2}}\Delta_{bh_{3}}c_{3}f(x)+\Delta_{-h_{1}}\Delta_{h_{2}}\Delta_{-ah_{3}}c_{4}f(y).$$
Letting $x-bh_{4}$ and $y-ah_{4}$ instead of $x$ and $y$ in
(2.1), we get
\begin{eqnarray*}
%
\Delta_{(a+b)h_{1}}\Delta_{(a-b)h_{2}}\Delta_{(2ab)h_{3}}c_{2}f(ax-by)&=&
\end{eqnarray*}
$$\Delta_{h_{1}}\Delta_{h_{2}}\Delta_{bh_{3}}c_{3}f(x-bh_{4})+\Delta_{-h_{1}}\Delta_{h_{2}}\Delta_{-ah_{4}}c_{4}f(y-ah_{4}).$$
Using this equality and (2.3), we find that
\begin{eqnarray}\label{6}
%
\Delta_{h_{1}}\Delta_{h_{2}}\Delta_{bh_{3}}\Delta_{-bh_{4}}c_{3}f(x)+\Delta_{-h_{1}}\Delta_{h_{2}}\Delta_{-ah_{3}}\Delta_{-ah_{4}}c_{4}f(y)
&=& 0
\end{eqnarray}
Letting $x+h_{5}$ and instead of $x$ and $y$ in (2.1), we get
$$\Delta_{h_{1}}\Delta_{h_{2}}\Delta_{bh_{3}}\Delta_{-bh_{4}}c_{3}f(x+h_{5})+\Delta_{-h_{1}}\Delta_{h_{2}}\Delta_{-ah_{3}}\Delta_{-ah_{4}}c_{4}f(y)=0.$$
finally Using this equality and (2.4) and since $c_{3}\neq0$,
then we obtain
$$\Delta_{h_{1}}\Delta_{h_{2}}\Delta_{bh_{3}}\Delta_{-bh_{4}}\Delta_{h_{5}}f(x)=0$$
for all $x\in V$. Therefore, by Theorem (1.1), $f$ is a GP
function of degree at most $4$. So there exist $a_{0}\in B$ and
symmetric $k$-additive functions $a_{k}:V^{k}\rightarrow B$ (for
$1\leq k\leq m$) such that
$$f(x)=a_{0}+\sum_{k=1}^{m}a_{k}^{\ast}(x)\ \ \ for\ all\ x\in V,$$
and $a^{\ast}_{m}\neq 0$, also
$$f(rx)=a_{0}+\sum_{k=1}^{m}r^{k}a_{k}^{\ast}(x)\ \ \ for\ all\ x\in V\ and\ r\in \mathbb{Q}.$$
The proof is complete.
\begin{thm}
Let $f:V\rightarrow B$ be a function such that $V$ and $B$ is a
vector space over $\mathbb{Q}$, $Df(x,y)=0$ for all $x, y\in V$ and
$f(0)=0$. If $|\frac{c_{3}+c_{5}+c_{6}}{c_{1}+c_{2}}|\neq0, 1$, then
$f$ is a monomial function of degree
$\log_{|a|}|\frac{c_{3}+c_{5}+c_{6}}{c_{1}+c_{2}}|$.
\end{thm}
\textbf{Proof.} by Theorem (2.1) $f$ is a GP function of
degree at most $4$ and there exist $a_{0}\in B$ and symmetric
$k$-additive functions $a_{k}:V^{k}\rightarrow B$ (for $1\leq k\leq
m$) such that
$$f(x)=a_{0}+\sum_{k=1}^{m}a_{k}^{\ast}(x)\ \ \ for\ all\ x\in V,$$
and $a^{\ast}_{m}\neq 0$, also
$$f(rx)=a_{0}+\sum_{k=1}^{m}r^{k}a_{k}^{\ast}(x)\ \ \ for\ all\ x\in V\ and\ r\in \mathbb{Q}.$$
Letting $y=0$ in (1.4), then
$$(c_{1}+c_{2})f(ax)=(c_{3}+c_{5}+c_{6})f(x)$$
for all $x\in V$. Thus,
$$f(ax)=\frac{c_{3}+c_{5}+c_{6}}{c_{1}+c_{2}} f(x)$$
for all $x\in V$. Since $|\frac{c_{3}+c_{5}+c_{6}}{c_{1}+c_{2}}|\neq
0, 1$, so the GP function $f$ must be a monomial function of degree
at most $4$. Therefore there exists a $k\in \mathbb{N}$ such that
$f=a_{k}^{\ast}$, where its a monomial function of degree $k$, so
$$|a|^{k}=|\frac{c_{3}+c_{5}+c_{6}}{c_{1}+c_{2}}|,$$
where it implies that $f$ is a monomial function of degree
$\log_{|a|}|\frac{c_{3}+c_{5}+c_{6}}{c_{1}+c_{2}}|$ and the proof is
complete.
\section{The generalized Hyers-Ulam stability of equation
(1.4)}\qquad

In the following, for the reader's convenience and explicit later
use, we will recall some fundamental results in fixed point theory.
\begin{defn}
The pair $(X, d)$ is called a generalized complete metric space if
$X$ is a nonempty set and $d:X^{2}\rightarrow [0,\infty]$ satisfies
the following conditions:
\begin{enumerate}
\item $d(x,y)\geq0$ and the equality holds if and only if $x=y$;
\item $d(x,y)=d(y,x)$;
\item $d(x,z)\leq d(x,y)+d(y,z)$;
\item every d-Cauchy sequence in X is d-convergent.
\end{enumerate}
for all $x, y\in X$.
\end{defn}
Note that the distance between two points in a generalized metric
space is permitted to be infinity.
\begin{defn}
Let $(X, d)$ be a metric space. A mapping $J:X\rightarrow X$
satisfies a Lipschitz condition with Lipschitz constant $L\geq0$ if
$$d(J(x),J(y))\leq Ld(x,y)$$
for all $x, y\in X$. If $L<1$, then $J$ is called a strictly
contractive map.
\end{defn}
\begin{thm}[38]\label{theor1}
Let $(X, d)$ be a generalized complete metric space and
$J:X\rightarrow X$ be strictly contractive mapping. Then for each
given element $x\in X$, either
$$d(J^{n}(x),J^{n+1}(x))=\infty$$
for all nonnegative integers $n$ or there exists a positive integer
$n_{0}$ such that
\begin{enumerate}
\item $d(J^{n}(x),J^{n+1}(x))<\infty$, for all $n\geq n_{0}$;
\item the sequence $\{J^{n}(x)\}$ converges to a fixed point $y^{\ast}$ of $J$;
\item $y^{\ast}$ is the unique fixed point of $J$ in the set $Y=\{y\in X\ :\ d(J^{n_{0}}(x),y)<\infty\}$;
\item $d(y,y^{\ast})\leq\frac{1}{1-L}d(J(y),y)$.
\end{enumerate}
\end{thm}
\begin{thm}\label{theo4}
Let $\lambda_{1}=\frac{c_{3}+c_{5}+c_{6}}{c_{1}+c_{2}}$ and
$a_{1}=a$, when $|\frac{c_{3}+c_{5}+c_{6}}{c_{1}+c_{2}}|>1$ and also
$\lambda_{2}=\frac{1}{\frac{c_{3}+c_{5}+c_{6}}{c_{1}+c_{2}}}$ and
$a_{2}=\frac{1}{a}$, when
$0<|\frac{c_{3}+c_{5}+c_{6}}{c_{1}+c_{2}}|<1$. Suppose that the
mapping $f:V\rightarrow X$ satisfies the conditions $f(0)=0$ and
\begin{eqnarray}\label{10}
  \|Df(x,y)\|&\leq& \psi(x,y)
\end{eqnarray}
for all $x,y\in V$, in which $\psi:V\times V\rightarrow [0,\infty)$
is a function. If there exists a positive real $L=L(i)<1$ $(i=1,2)$
such that
\begin{eqnarray}\label{7}
%
\lim_{n\rightarrow\infty} \lambda_{i}^{-n}\psi(a_{i}^{n}x,a_{i}^{n}y)&=&0,\\
\psi(a_{i}x,0)&\leq&L_{i}|\lambda_{i}|\psi(x,0)\label{9}
\end{eqnarray}
for all $x,y\in V$, then there exists a unique $T:V\rightarrow X$
such that $DT(x,y)=0$ and
$$\|f(x)-T(x)\|\leq\frac{L^{2-i}}{(1-L)|c_{3}+c_{5}+c_{6}|}\psi(x,0)$$
for all for $x,y\in V$.
\end{thm}
\textbf{Proof.} Let us consider the set $A:=\{g:V\rightarrow X\}$
and introduce the generalized metric on $A$:
$$d(g,h)=\sup_{\{x\in X\ :\ \psi(x,0)\neq 0\}}\frac{\|g(x)-h(x)\|}{\psi(x,0)}.$$
It is easy to show that $(A, d)$ is generalized complete metric
space. Now we define the function $J:A\rightarrow A$ with
\begin{eqnarray}
%
J(g(x))&=& \frac{1}{\lambda_{i}}g(a_{i}x)
\end{eqnarray}
for all $g\in A$ and $x\in V$. We can write,
\begin{eqnarray*}
%
d(J(g),J(h))&=&\sup_{\{x\ :\ \psi(x,0)\neq0\}}\frac{\|g(a_{i}x)-h(a_{i}x)\|}{|\lambda_{i}|\psi(x,0)} \\
&\leq& L\sup_{\{x\ :\
\psi(x,0)\neq0\}}\frac{\|g(a_{i}x)-h(a_{i}x)\|}{\psi(a_{i}x,0)}=Ld(g,h),
\end{eqnarray*}
that is $J$ is a strictly contractive selfmapping of $A$, with the
Lipschitz constant $L$. We set $y=0$ in the hypothesis (3.1),
then we obtain
\begin{eqnarray}\label{8}
%
\|\frac{f(ax)}{\lambda}-f(x)\| &\leq&
\frac{1}{|c_{3}+c_{5}+c_{6}|}\psi(x,0)
\end{eqnarray}
for all $x\in V$. Now if $i=1$, then from (3.5), we get
$$\|\frac{f(a_{1}x)}{\lambda_{1}}-f(x)\|\leq\frac{1}{|c_{3}+c_{5}+c_{6}|}\psi(x,0)$$
for all $x\in V$ and it implies that
$d(J(f),f)<\frac{1}{|c_{3}+c_{5}+c_{6}|}<\infty$. If $i=2$, then
from (3.5) and (3.3) , we get
$$\|\frac{f(a_{2}x)}{\lambda_{2}}-f(x)\|\leq\frac{L}{|c_{3}+c_{5}+c_{6}|}\psi(x,0)$$
for all $x\in V$ and it implies that
$d(J(f),f)<\frac{L}{|c_{3}+c_{5}+c_{6}|}<\infty$. By Theorem
(3.3), there exists a mapping $T:V\rightarrow X$ such that
\begin{enumerate}
\item $T$ is a fixed point of $J$, i.e.,
\begin{eqnarray}
%
T(a_{i}x) &=& \lambda_{i}T(x)
\end{eqnarray}
for all $x\in V$. The mapping $T$ is a unique fixed point of $J$ in
the set $\tilde{A}=\{h\in A\ :\ d(f,h)<\infty\}$.
\item $d(J^{n}(f),T)\rightarrow 0$ as $n\rightarrow \infty$. This implies that
$$T(x)=\lim_{n\rightarrow\infty}\frac{f(a_{i}^{n}x)}{\lambda_{i}^{n}}$$
for all $x\in V$ and also
$$\lim_{n\rightarrow\infty}\frac{Df(a_{i}^{n}x,a_{i}^{n}y)}{\lambda_{i}^{n}}=DT(x,y)$$
\item $d(f,T)\leq\frac{1}{1-L}d(J(f),f)$ and $d(J(f),f)\leq \frac{L^{2-i}}{|c_{3}+c_{5}+c_{6}|}$, so
$$d(f,T)\leq\frac{L^{1-i}}{(1-L)|c_{3}+c_{5}+c_{6}|}.$$
\end{enumerate}
Now from (3.1), we see that
\begin{eqnarray}
%
\|\frac{Df(a_{i}^{n}x,a_{i}^{n}y)}{\lambda_{i}^{n}}\|&\leq&
\frac{\psi(a_{i}^{n}x,a_{i}^{n}y)}{|\lambda_{i}|^{n}}
\end{eqnarray}
for all $x,y\in V$ and letting $n$ to infinity, we get $DT(x,y)=0$
for all $x,y\in V$ and the proof is complete.
\begin{thm}
Let $\gamma_{1}=\frac{c_{3}+c_{4}+c_{6}}{c_{1}+c_{2}}$ and
$k_{1}=a+b$, when $|\frac{c_{3}+c_{4}+c_{6}}{c_{1}+c_{2}}|>1$ and
also $\gamma_{2}=\frac{1}{\frac{c_{3}+c_{4}+c_{6}}{c_{1}+c_{2}}}$
and $k_{2}=\frac{1}{a+b}$, when $0<|\gamma|<1$. Suppose that the
mapping $f:V\rightarrow X$ satisfies the conditions $f(0)=0$ and
\begin{eqnarray}\label{c10}
%
\|Df(x,y)\|&\leq& \psi(x,y)
\end{eqnarray}
for all $x,y\in V$ and also $c_{2}=0$ in the (1.4), in which
$\psi:V\times V\rightarrow [0,\infty)$ is a function. If there
exists a positive real $L=L(i)<1$ $(i=1,2)$ such that
\begin{eqnarray}\label{c7}
%
\lim_{n\rightarrow\infty} \gamma_{i}^{-n}\psi(k_{i}^{n}x,k_{i}^{n}y)&=&0,\\
\psi(k_{i}x,k_{i}x)&\leq&L_{i}|\gamma_{i}|\psi(x,x)\label{c9}
\end{eqnarray}
for all $x,y\in V$, then there exists a unique $T:V\rightarrow X$
such that $DT(x,y)=0$ and
$$\|f(x)-T(x)\|\leq\frac{L^{2-i}}{(1-L)|c_{3}+c_{4}+c_{6}|}\psi(x,x)$$
for all for $x,y\in V$.
\end{thm}
\textbf{Proof.} Let us consider the set $U:=\{g:V\rightarrow X\}$
and introduce the generalized metric on $U$:
$$d(g,h)=\sup_{\{x\in X\ :\ \psi(x,x)\neq 0\}}\frac{\|g(x)-h(x)\|}{\psi(x,x)}.$$
It is easy to show that $(U, d)$ is generalized complete metric
space. Now we define the function $J:A\rightarrow A$ with
\begin{eqnarray}
%
J(g(x))&=& \frac{1}{\gamma_{i}}g(k_{i}x)
\end{eqnarray}
for all $g\in U$ and $x\in V$. We can write,
\begin{eqnarray*}
%
d(J(g),J(h))&=&\sup_{\{x\ :\ \psi(x,x)\neq0\}}\frac{\|g(k_{i}x)-h(k_{i}x)\|}{|\gamma_{i}|\psi(x,x)} \\
&\leq& L\sup_{\{x\ :\
\psi(x,x)\neq0\}}\frac{\|g(k_{i}x)-h(k_{i}x)\|}{\psi(k_{i}x,k_{i}x)}=Ld(g,h),
\end{eqnarray*}
that is $J$ is a strictly contractive selfmapping of $U$, with the
Lipschitz constant $L$. We set $y=x$ in the hypothesis (3.8),
then we obtain
\begin{eqnarray}\label{c8}
%
\|\frac{f((a+b)x)}{\frac{c_{3}+c_{4}+c_{6}}{c_{1}+c_{2}}}-f(x)\|
&\leq& \frac{1}{|c_{3}+c_{4}+c_{6}|}\psi(x,x)
\end{eqnarray}
for all $x\in V$. Now if $i=1$, then from (3.12), we get
$$\|\frac{f(k_{1}x)}{\gamma_{1}}-f(x)\|\leq\frac{1}{|c_{3}+c_{4}+c_{6}|}\psi(x,x)$$
for all $x\in V$ and it implies that
$d(J(f),f)<\frac{1}{|c_{3}+c_{4}+c_{6}|}<\infty$. If $i=2$, then
from (3.12) and (3.10) , we get
$$\|\frac{f(k_{2}x)}{\gamma_{2}}-f(x)\|\leq\frac{L}{|c_{3}+c_{4}+c_{6}|}\psi(x,x)$$
for all $x\in V$ and it implies that
$d(J(f),f)<\frac{L}{|c_{3}+c_{4}+c_{6}|}<\infty$. By Theorem
(3.3), there exists a mapping $T:V\rightarrow X$ such that
\begin{enumerate}
\item $T$ is a fixed point of $J$, i.e.,
\begin{eqnarray}
%
T(k_{i}x) &=& \gamma_{i}T(x)
\end{eqnarray}
for all $x\in V$. The mapping $T$ is a unique fixed point of $J$ in
the set $\tilde{U}=\{h\in U\ :\ d(f,h)<\infty\}$.
\item $d(J^{n}(f),T)\rightarrow 0$ as $n\rightarrow \infty$. This implies that
$$T(x)=\lim_{n\rightarrow\infty}\frac{f(k_{i}^{n}x)}{\gamma_{i}^{n}}$$
for all $x\in V$ and it implies that
$$\lim_{n\rightarrow\infty}\frac{Df(k_{i}^{n}x,k_{i}^{n}y)}{\gamma_{i}^{n}}=DT(x,y)$$
\item $d(f,T)\leq\frac{1}{1-L}d(J(f),f)$ and $d(J(f),f)\leq \frac{L^{1-i}}{|c_{3}+c_{4}+c_{6}|}$ , which implies,
$$d(f,T)\leq\frac{L^{2-i}}{(1-L)|c_{3}+c_{4}+c_{6}|}.$$
\end{enumerate}
Now from (3.1), we get
\begin{eqnarray}
%
\|\frac{Df(k_{i}^{n}x,k_{i}^{n}y)}{\gamma_{i}^{n}}\|&\leq&
\frac{\psi(k_{i}^{n}x,k_{i}^{n}y)}{|\gamma_{i}|^{n}}
\end{eqnarray}
for all $x,y\in V$ and letting $n$ to infinity, we get $DT(x,y)=0$
for all $x,y\in V$ and the proof is complete.
\bibliographystyle{amsplain}

\end{document}